\documentclass{amsart}
\theoremstyle{plain}
\newtheorem{thm}{Theorem}
 
 \newtheorem{definition}[thm]{Definition}

\errorcontextlines=0 \numberwithin{equation}{section}

\begin{document}
\large
\title[Ricci flow with hyperbolic metrics]
{Ricci Flow with hyperbolic warped product metrics}
\author{Li Ma and  Xingwang Xu}

\address{Ma: Department of Mathematical Sciences \\
Tsinghua University \\
Beijing 100084 \\
China}

\email{lma@math.tsinghua.edu.cn}

\address{Xu: Mathematics Department\\
The National University of Singapore\\
10 Kent Ridge Crescent\\
Singapore 119260}

\email{matxuxw@math.nus.edu.sg}

\dedicatory{}
\date{}
\thanks{The research is partially supported by the National Natural Science
Foundation of China 10631020 and SRFDP 20060003002}


\begin{abstract}
In this short note, we show that the negative curvature is
preserved in the deformation of hyperbolic warped product metrics
under Ricci flow. It is also showed that the flow converges to a
flat metric as time going to infinity.
\end{abstract}
\maketitle

\section{Introduction}

In this short note, we show by an example that the negative
curvature is preserved in the deformation of hyperbolic warped
product metrics under Ricci flow. It is also showed that the Ricci
flow converges to a flat metric as time going to infinity. Given a
time interval $I=(a,b)$. Recall that the Ricci flow on a
non-compact manifold $M^{n+1}$ is defined by a one parameter
family of metrics $\{g(t); t\in I\}$ satisfying
\begin{equation}\label{Rcflow}
\partial_tg(t)=-2Rc(g(t)),\; \; t\in I
\end{equation}
where $Rc(g(t))$ is the Ricci tensor of the metric $g(t)$. We
denote by $R(t)$ the scalar curvature of $g(t)$. Then the scalar
curvature function $R(t)$ satisfies the following evolution
equation
\begin{equation}\label{Scflow}
\partial_t R=\Delta R+2|Rc|^2\geq \Delta R+\frac{2}{n}R^2.
\end{equation}
One often use this evolution equation to get the finite blow up
for Ricci flow the positive curvature. We shall consider the flow
of warped product metrics and use this equation to get a lower
bound which forces the metrics $g(t)$ going to a flat metric. In
many physical important examples of Einstein metrics, we meet the
warped product metrics. So it is nature to consider the Ricci flow
for warped product metrics.

 Let
$M=\mathbb{R}_+\times N^n$ with the product metric
$$
g=\phi(x)^2dx^2+\psi(x)^2\hat{g}.
$$
Here $(N^n,\hat{g})$ is an Einstein manifold of dimension $n\geq
2$, $\phi(x)$ and $\psi(x)$ are two smooth positive functions of
the variable $x>0$. In this paper, we take $N=S^n$ be the standard
n-sphere so that $\hat{g}=d\sigma^2$. If $n=1$, the related
results on Ricci flow can be found in \cite{M}.

Introduce the arc-length parameter $s$ and sectional curvatures in
the following way. Let
\begin{equation}\label{length}
s=\int_0^x\phi(\tau)d\tau.
\end{equation}
Let $$ K_0=-\frac{\psi_{ss}}{\psi}.
$$ be the sectional curvature of the 2-planes perpendicular to the
spheres $\{x\}\times S^n$
 and
$$
K_1=\frac{1-\psi_s^2}{\psi^2}
$$
the sectional curvature of the 2-planes tangential to the spheres.
Note that
$$
\frac{\partial}{\partial
s}=\frac{1}{\phi(x)}\frac{\partial}{\partial x}.
$$

In terms of the arc-length parameter $s$, the metric can be read
as
$$ g=ds^2+\psi^2\hat{g}.
$$

Recall that the standard hyperbolic n-space can be written as
$$
(\mathbb{R_+}\times S^{n-1}, dx^2+\sinh^2x \hat{g}).
$$
Using this model, we give the following definition.

\begin{definition} We say the product metric $g$ on the manifold $M$ is of hyperbolic
type if it has bounded negative sectional curvature with
$$
\frac{\phi(x)}{x}\approx 1, \; \;
$$
and $$ \frac{\psi(x)}{\sinh(x)}\approx 1,\;
$$ as $x\to +\infty$, where $f\approx 1$ means that there exist two positive constants
$C_1$ and $C_2$ such that $C_1\leq f\leq C_2$ uniformly for $x$
large.
 We write by $HM(n+1)$ the class of product
metrics of hyperbolic type on $M$.
\end{definition}
We remark that if we impose the conditions
$$
\partial_x^{\alpha}({\phi(x)}/x)\to 0, \; \;
\partial_x^{\alpha}(\frac{\psi(x)}{\sinh(x)})\to 0,
$$
for $x\to +\infty$, where $1\leq \alpha\leq 2$, then we can have
that
$$
K_0\to -1, \; \; K_1\to -1
$$
as $|x|\to \infty$.

 It is known
that the Ricci curvature tensor of the product metric $g$ is given
by
$$
Rc=n \{-\frac{\psi_{xx}\psi+\psi_x\phi_x}{\psi\phi}
\}dx^2+\{-\frac{\phi\psi_x\psi_{xx}-(n-1)\phi\psi_x^2+\psi\phi_x\psi_x}{\phi^3}+n-1\}\hat{g}.
$$

 Let
$$ g_s=\psi^2\hat{g}.
$$
Then we have the hessian formula for $s$ (see \cite{PP}):
$$
Hess (s)=\psi \psi_s \hat{g}=\frac{\psi_s}{\psi}g_s,
$$
and the Ricci tensor can be written as
$$
Rc=nK_0ds^2+[K_0+(n-1)K_1]\psi^2\hat{g}.
$$
See \cite{SM} and \cite{AK}.

Note that the scalar curvature is given by
$$
R=nK_0+n[K_0+(n-1)K_1].
$$

Then for the warped product metrics, the Ricci flow (\ref{Rcflow})
becomes the system
\begin{equation}\label{Rcst}
\left\{\begin{array}{ll}
\psi_t=-[K_0+(n-1)K_1]\psi=\psi_{ss}-(n-1)\frac{1-\psi_s^2}{\psi},
 \\
\phi_t =-nK_0\phi=n\frac{\psi_{ss}}{\psi}\phi.
\end{array}
\right.
\end{equation}
In the xt coordinates, we have
\begin{equation}\label{Rcxt}
\left\{\begin{array}{ll} \psi_t=
\frac{1}{\phi(x)}\frac{\partial}{\partial
x}\frac{1}{\phi(x)}\frac{\partial}{\partial x}\psi
-(n-1)\frac{1-(\frac{1}{\phi(x)}\frac{\partial}{\partial
x}\psi)^2}{\psi},
 \\
\phi_t = n\frac{\frac{\partial}{\partial
x}\frac{1}{\phi(x)}\frac{\partial}{\partial x}\psi}{\psi}.
\end{array}
\right.
\end{equation}

We prove the following:
\begin{thm} \label{thm1} Let $n\geq 2$. Given a complete initial metric $g_0\in HM(n+1)$.
Then the Ricci flow exists globally and keeps in the class
$HM(n+1)$; furthermore, as $t\to \infty$, the flow converges
uniformly to the flat metric.
\end{thm}
We note that the local existence of Ricci flow has been proved by
W.Shi \cite{Shi89} and the uniqueness result can be proved in the
same way as in Theorem 3.3 in M.Simon \cite{SM} by using the
harmonic map heat flow associated to the Ricci flow (see
R.Hamilton \cite{H95}). We remark that in our case, the condition
(4) in \cite{SM} has not been satisfied. However, the argument
there can be modified to our hyperbolic case. The goal in the
papers \cite{SM} and \cite{AK} (see also \cite{CK}) is to show
that the pinching singularity in finite time along the Ricci flow
occurs for some class of warped product metrics with non-negative
curvature, as one can expect from the intuition. Although it is
conjectured by R.Hamilton that the globally Ricci flow exists for
the negative sectional curvature case, the convergence of the
Ricci flow can not be easily classified.

We now give some remark about the possible application of our
result to study the rigidity problem of hyperbolic warped product
metrics. It may also be useful for our result to study the
evolution of the quasi-local mass of Brown and York on the spheres
(when $n=2$). To introduce the Brown-York mass of geodesic
spheres, we look at the evolution of their mean curvatures.

 For a fixed $r>0$,
let
$$
\Sigma:=\partial B(r)=\{p\in M; s(p)=r\}
$$
be the sphere of radius $r$ in the warped Riemannian manifold
$(M,g)$. Note that the volume of $\partial B(r)$ is
$$
V(r)=|\Sigma|=\int_{\Sigma}d\sigma=\int_N\psi(s)^{n}dv_{\hat{g}}=\psi(r)^nV(N).
$$
Here $V(N)$ is the volume of $N=S^n$ in the metric $\hat{g}$.

 By the hessian formula for the function $s$, we have the mean
curvature of the sphere $\Sigma$ given by
\begin{equation}\label{hessian}
H=\Delta s=n\frac{\psi_s}{\psi}=\frac{V'(s)}{V(s)}.
\end{equation}

Let $n=2$. Let's recall the definition of the Brown-York
quasi-local mass (see \cite{B93} for definition and see \cite{W02}
for more)
$$
m(r):=m(\Sigma)=\int_{\Sigma}(H_0-H).
$$
where $H_0$ is the mean curvature of $\Sigma$ in the standard
hyperbolic space. Recall that
$$ H_0(x)=\frac{n\cosh(x)}{\sinh(x)}.
$$
In our case, it is clear that
$$
m(r):=m(\Sigma)=(H_0-H)V(s)=H_0V(r)-V_s(r).
$$

Motivated by results in \cite{P}, \cite{CM05}, and \cite{DM1}, it
is possible to prove some (weighted) monotonicity formula for the
mass $m(r)$ along the Ricci flow.

\section{Basic evolution equations}
In this section, we derive the evolution equations for the
arc-length parameter $s=s(x,t)$, the volume function $V$, the
curvature function $K_0$, and mean curvature function $H$ along
the Ricci flow (\ref{Rcst}). Some of these evolution equations for
the derivatives of $\psi$ were derived in \cite{AK}. We derive
more formulae than what we can use in this paper since they may be
useful to study the evolutions of quasi-local masses and ADM mass
(see \cite{B93}, \cite{DM1}, and \cite{W02} for the definition)
along the Ricci flow.

Recall that the
 product metric $g(t)$ is given by
 $$
g(t)=\phi(x,t)^2dx^2+\psi(x,t)^2\hat{g},
 $$
 which is
 $$
g(t)=ds^2+\psi(s,t)^2\hat{g}
 $$
 in the st-coordinates. For a given smooth radial function
 $f=f(s)$, we have
 $$
\Delta
f=\Delta_{g(t)}f=\frac{1}{\phi^2}f_{xx}+\frac{1}{\phi\psi^n}(\frac{\psi^n}{\psi})_xf_x
 $$
 or
 $$
\Delta f=f_{ss}+n\frac{\psi_s}{\psi}f_s.
 $$

First of all, along the Ricci flow, we have the following
commutator relation (see \cite{MC}):
\begin{equation}\label{commutator}
[\frac{\partial}{\partial t},\frac{\partial}{\partial
s}]=\frac{\partial}{\partial t}\frac{\partial}{\partial
s}-\frac{\partial}{\partial s}\frac{\partial}{\partial
t}=-\frac{1}{\phi}\frac{\partial \phi}{\partial
t}\frac{\partial}{\partial
s}=-n\frac{\psi_{ss}}{\psi}\frac{\partial}{\partial
s}=nK_0\frac{\partial}{\partial s}.
\end{equation}

We then consider the evolution equation for the arc-length
parameters of the flow. By the equation (\ref{length}), we have
$$
\frac{\partial s}{\partial t}=\int^x
\phi_td\xi=n\int^s\frac{\psi_{ss}}{\psi}d\tau.
$$
By this, we have
$$
\frac{\partial s}{\partial
t}=n\frac{\psi_{s}}{\psi}+n\int^s\frac{\psi_{s}^2}{\psi^2}.
$$
Using the expression (\ref{hessian}) for the mean curvature $H$,
we have
\begin{equation}\label{sderivative}
\frac{\partial s}{\partial t}=H+\frac{1}{n}\int^s H^2ds=\Delta
s+\frac{1}{n}\int^s H^2d\tau.
\end{equation}

 Note that under the Ricci flow (\ref{Rcst}), the volume
$V(t)$ is changed by
$$
V_t=n\int_N \psi^{n-1}\psi_tdv_{\hat{g}}=nV(N)
\psi^{n-1}(\psi_{ss}-(n-1)\frac{1-\psi_s^2}{\psi}).
$$
Then we can simplify it to the following form:
\begin{equation}\label{Vflow}
V_t=V_{ss}-\frac{n-1}{\psi^2}V=V_{ss}-(n-1)V(N)^{2/n}V^{(n-2)/n},
\end{equation}
which is a parabolic equation for the volume $V(s,t)$ of the
spheres $\Sigma$. This gives us that
$$
(\partial_t-\Delta
)V=-\frac{n\psi_s}{\psi}V_s-(n-1)V(N)^{2/n}V^{(n-2)/n}.
$$
Note that, when $n=2$, we have
$$
V_t=V_{ss}-V(N).
$$

Using the Ricci flow equation (\ref{Rcst}), we have
$$
(\partial_t-\Delta)\psi=K_1\psi-\frac{n}{\psi}.
$$

Using the commutator relation (\ref{commutator}), we have
$$
\partial_t\psi_s=\partial_s\psi_t-\frac{1}{\phi}\phi_t\psi_s.
$$
Then we have
\begin{equation}\label{first}
\partial_t\psi_s-(\psi_s)_{ss}=(n-2)\frac{\psi_{ss}\psi_s}{\psi}+(n-1)\frac{1-\psi_s^2}{\psi^2}\psi_s,
\end{equation}
which can also be written as
$$
\partial_t\psi_s-(\psi_s)_{ss}=[-(n-2)K_0+(n-1)K_1]\psi_s.
$$
Here we have used that
$$
(K_0)_s=-\frac{\psi_{sss}}{\psi}-K_0\frac{\psi_s}{\psi}, \; \;
(K_1)_s=-\frac{\psi_s}{\psi}(K_1-K_0).
$$
Hence, by (\ref{first}), we have
\begin{equation}\label{first2}
\partial_t\psi_s-\Delta(\psi_s)=[2K_0+(n-1)K_1]\psi_s.
\end{equation}

 Let $w=\psi_{ss}$. Then we have
\begin{equation}\label{second}
w_t-w_{ss}=(n-2)\frac{\psi_s}{\psi}w_s
-(2K_0-(4n-5)\frac{\psi_s^2}{\psi^2}+\frac{n-1}{\psi^2})w-2(n-1)K_1\frac{\psi_s^2}{\psi}.
\end{equation}
Hence,
$$
(\partial_t-\Delta)w=-2\frac{n\psi_s}{\psi}w_s
-(2K_0-(4n-5)\frac{\psi_s^2}{\psi^2}+\frac{n-1}{\psi^2})w-2(n-1)K_1\frac{\psi_s^2}{\psi}.
$$

 We now compute the evolution equation for the
curvature function $$K=-K_0=\psi_{ss}/\psi,$$
which is more
important to us. Using (\ref{second}), we can derive that
\begin{equation}\label{Kflow}
K_t-K_{ss}=
n\frac{n\psi_s}{\psi}K_s-2K^2-4(n-1)K_1K-2(n-1)K_1\frac{\psi_s^2}{\psi^2}.
\end{equation}
That is to say,
$$
K_t-\Delta K= -4K^2-2(n-1)K_1K-2(n-1)K_1\frac{\psi_s^2}{\psi^2}.
$$

We note that
$$
H_t=-n\psi_s\frac{\psi_t}{\psi^2}+n\frac{\partial_t\partial_s\psi}{\psi}=
-n\psi_s\frac{\psi_t}{\psi^2}+n\frac{1}{\psi}(\partial_s\psi_t-\frac{1}{\phi}
\phi_t\psi_s).
$$
Using the relations that
$$
\frac{\psi_{ss}}{\psi}=\frac{H_s}{n}+\frac{H^2}{n^2}
$$
and
$$
\frac{\psi_{sss}}{\psi}=\frac{H_{ss}}{n}+\frac{3\psi_s\psi_{ss}}{\psi^2}-\frac{2H^3}{n^3}
=\frac{H_{ss}}{n}+\frac{3HH_s}{n^2}+\frac{H^3}{n^3},
$$
we can easily find that
\begin{equation}\label{Hflow}
H_t=H_{ss}+\frac{2n-1}{n}HH_s-\frac{1}{n^2}H^3+2(n-1)V(N)^{2/n}V^{-2/n}H.
\end{equation}
Then we have
$$
(\partial_t-\Delta)H=\frac{n-1}{n}HH_s-\frac{1}{n^2}H^3+2(n-1)V(N)^{2/n}V^{-2/n}H.
$$

Let
$$
a(x,t)=\psi^2(K_1-K_0)
$$
be the measure of the difference of the two kinds of sectional
curvatures. Then we have $$
a_t=a_{ss}+(n-4)\frac{\psi_s}{\psi}a_s-4(n-1)\frac{\psi_s^2}{\psi^2}a,
$$
which is
$$
(\partial_t-\Delta)a=-4\frac{\psi_s}{\psi}a_s-4(n-1)\frac{\psi_s^2}{\psi^2}a.
$$
 We can use the maximum principle to this equation to find a pinching estimate.

\section{Proof of Theorem \ref{thm1}}
 In this section we mainly prove Theorem \ref{thm1}.

For any given smooth invertible function
$f:\mathbb{R}\to\mathbb{R}$, we consider the evolution equations
for the function $u:=f(\psi)$. Note that
$$
u_t=f'\psi_t, \; \; u_s=f'\psi_s,
$$
and
$$
u_{ss}=f'\psi_{ss}+f^{''}\psi_s^2.
$$
Then
$$
u_t-u_{ss}=f'(\psi_t-\psi_{ss})-\frac{f^{''}\circ f^{-1}}{|f'\circ
f^{-1}|^2}u_s^2=-f'((n-1)K_1+n)\psi-\frac{f^{''}\circ
f^{-1}}{|f'\circ f^{-1}|^2}u_s^2,
$$
and
\begin{equation}\label{Bflow}
u_t-u_{ss}=-((n-1)K_1+n)f'\circ
f^{-1}(u)f^{-1}(u)-\frac{f^{''}\circ f^{-1}}{|f'\circ
f^{-1}|^2}u_s^2.
\end{equation}
We shall use this evolution equation to study the behavior of the
function $\psi$. In particular, for $u=f(u)=\psi^k$ for $\psi<1$
and $u=1$ for $\psi\geq 1$, we have that $f^{-1}=u^{1/k}$,
$f'=k\psi^{k-1}$, $f^{''}=k(k-1)\psi^{k-2}$, and
$$
u_t-u_{ss}=-k((n-1)K_0+n)u-\frac{k-1}{k}u^{-1}u_s^2.
$$
If $k$ is a negative integer, we know that from the argument as in
\cite{DM1} that $\psi(x,t)$ has the same decay as $\psi(x,0)$ at
infinity.

We now need to study the behavior of functions $\psi$ and $\phi$
at $x=0$. From the definition of Riemannian metric, we know that
for all $t>0$,
$$
s(0,t)=0, \; s_x(0,t)=0=\phi(0,t),\; \phi_x(0,t)=1 \; \psi(0,t)=0,
\; \psi_s(0,t)=1.
$$
Using the formulae for curvature $K_0$ and bounded curvature
conditions we have
$$
\psi_{ss}(0,t)=0.
$$

Using the maximum principle (see \cite{Shi89}) to the evolution
equation (\ref{first}), we know that
$$
\psi_s>1, \; \; for \; \; x>0, \; t>0.
$$
This implies that $K_1<0$ for $x>0$ and $t>0$.

Using the maximum principle to the evolution equation
(\ref{second}), we know that
$$
\psi_{ss}>0, \; \; for \; \; x>0, \; t>0,
$$
which gives us that $K_0<0$ for $x>0$ and $t>0$. Then we have
$R(t)<0$ for all $t>0$.

Taking s-derivative in the formulae $\psi_{ss}=-K_0\psi$, we get
that $\psi_{sss}(0,t)=-K_0(0,t)\psi_s(0,t)\geq 0 $.

 Using the maximum principle to the evolution equation
(\ref{Rcxt}), we know that
$$
\psi>0,\; \;  \phi>0
$$
for all $t>0$ and $x>0$.

Using the negativity of the curvatures $K_0$ and $K_1$, the scalar
curvature evolution equation (\ref{Scflow}), and the maximum
principle, we have the following uniform lower bound of the scalar
curvature $R$,
\begin{equation}\label{scalarbound}
0> R(t)\geq \frac{-1}{2t/n-(\inf_M R(0))^{-1}}.
\end{equation}
Hence we have the uniform bound for the sectional curvature
functions $K_0$ and $K_1$. Using Shi's estimate (see Theorem 1.1
in \cite{Shi89}), we know that the flow $g(t)$ exists globally,
and by our estimate (\ref{scalarbound}), as $t\to \infty$, the
flow $g(t)$ converges to a smooth metric $g_{\infty}$ with the
zero sectional curvature.

We now consider the behavior of the solution pair $(\phi, \psi))$
to the flow (\ref{Rcst}). By (\ref{Rcst}), we have that
$$
\partial_t(\log \psi)=-[K_0+(n-1)K_1], \; \; \partial_t(\log
\phi)=-nK_0.
$$
Integrating over the time interval $[0,t]$, we get
$$
\psi(x,t)=\psi(x,t)exp(-\int_0^t[K_0+(n-1)K_1]),
$$
and $$
 \phi(x,t)=\phi(x,0)exp(-\int_0^tK_0),
$$
which say that the behavior of $\psi(x,t)$ and $\phi(x,t)$ are
equivalent to $\psi(x,0)$ and $\phi(x,0)$ respectively at infinity
$|x|\to \infty$.

This proves that the Ricci flow keeps the class $HM(n+1)$.

 By this, we complete the proof of Theorem
\ref{thm1}.

\section{Discussions}

Looking at the example of the flow $(g(t))$ defined by
$$
\partial_tg(t)=2g(t).
$$
Then we have
$$
g(t)=e^{2t}g(0).
$$
Hence, as $t\to +\infty$, we have
$$
g(t)\to g_{\infty},
$$
where $g_{\infty}$ is a flat metric. In this sense, Theorem
\ref{thm1} is not surprising.

 To make the limit metric of the Ricci flow to have negative
curvature, we may study the following modified Ricci flow:
\begin{equation}\label{mRcflow}
\partial_tg(t)=-2Rc(g(t))-2ng,\; \; t\in I
\end{equation}

If $g(t)$ is a modified Ricci flow, then the family defined by $$
\bar{g}(\tau) =(1+n\tau)g((\frac{1}{n}\log(1+n\tau)))
$$
is a Ricci flow. Hence, we can use our Theorem \ref{thm1} to get
globally existence of the modified Ricci flow (\ref{mRcflow}).
However, we shall leave the convergence topic of the global
modified Ricci flow to future study.

{\bf Acknowledgement}. The paper has been done when the first
named author is visiting the Department of Mathematics  in NUS,
Singapore. The author thanks the Department of Mathematics  in
NUS, Singapore for hospitality.

\end{document}